\documentclass[10pt]{article}
\usepackage{amssymb, amscd, verbatim, amsthm}
\usepackage[arrow, matrix, curve]{xy}
\usepackage{hyperref}
\usepackage{amsmath}
\usepackage{geometry}
\usepackage{color}
\usepackage{mathrsfs}
\usepackage{appendix}

\newtheorem{leer}{}[section]
\newtheorem{thm}[leer]{Theorem}

\newtheorem{defi}[leer]{Definition}
\newtheorem{rema}[leer]{Remark} 
\newtheorem{prop}[leer]{Proposition}
\newtheorem{lemm}[leer]{Lemma}
\newtheorem{coro}[leer]{Corollary}

\newtheorem*{mthm}{Main Theorem}

\newcommand{\plim}{\mathop{\varprojlim}\limits}

\newcommand{\OO}{{\cal O}}

\newcommand{\Gal}{{\mathrm{Gal}}}
\newcommand{\Spec}{\mathrm{Spec}}
\newcommand{\GL}{\mathrm{GL}}
\newcommand{\Hom}{\mathrm{Hom}}
\newcommand{\trdeg}{\mathrm{trdeg}}
\newcommand{\mf}{\mathfrak}
\newcommand{\ol}[1]{\overline{#1}}
\newcommand{\ul}[1]{\underline{#1}}
\newcommand{\Aut}{\mathrm{Aut}}
\newcommand{\End}{\mathrm{End}}

\newcommand{\im}{\mathrm{im}}
\newcommand{\chara}{\mathrm{char}}

\newcommand{\orho}{\ol{\rho}}

\newcommand{\Lie}{\mf{Lie}}
\newcommand{\ab}{\mathrm{ab}}
\newcommand{\solv}{\mathrm{solv}}
\newcommand{\cyc}{\mathrm{cyc}}
\newcommand{\FSQ}{\mathrm{FSQ}}
\newcommand{\JH}{\mathrm{JH}}

\def\Pp{{\mathbb{P}}}
\def\Ff{{\mathbb{F}}}

\def\Qq{{\mathbb{Q}}}
\def\Zz{{\mathbb{Z}}}

\def\Nn{{\mathbb{N}}}
\def\Gg{{\mathbb{G}}}
\def\Ll{{\mathbb{L}}}

\begin{document}

\baselineskip=17pt

\title{Group theoretical independence of $\ell$-adic Galois representations}
\author{Sebastian Petersen}

\parindent0em
\parskip1em


\maketitle

\renewcommand{\thefootnote}{}



\renewcommand{\thefootnote}{\arabic{footnote}}
\setcounter{footnote}{0}


\begin{abstract}

Let $K/\Qq$ be a finitely generated field of characteristic zero and 
$X/K$ a smooth projective variety. Fix $q\in\Nn$. For every prime number $\ell$ 
let $\rho_\ell$ be the representation of $\Gal(K)$ on the \'etale cohomology group
$H^q(X_{\ol{K}}, \Qq_\ell)$. For a field $k$ we denote by $k_{\ab}$ its maximal abelian Galois 
extension. We prove that 
there exist finite Galois extensions $k/\Qq$ and $F/K$ such that 
the restricted family of representations $(\rho_\ell|\Gal(k_{\ab} F))_\ell$ is group theoretically
independent in the sense that $\rho_{\ell_1}(\Gal(k_{\ab} F))$ and $\rho_{\ell_2}(\Gal(k_{\ab} F))$ 
do not have a common finite simple quotient group for all prime numbers $\ell_1\neq \ell_2$. 
\end{abstract}


\section*{Introduction} 
Following Serre (cf. \cite{bible}), let us call an arbitrary family $(\rho_\ell: G\to G_\ell)_{\ell\in\Ll}$ of continuous homomorphisms of 
profinite groups
{\em independent}, if the induced homomorphism
$\rho: G\to \prod_{\ell\in\Ll} \rho_\ell(G)$ is surjective. The family $(\rho_\ell: G\to G_\ell)_{\ell\in\Ll}$
is said to be {\em group theoretically independent} if for all prime numbers $\ell_1\neq \ell_2$ the
groups $\rho_{\ell_1}(\Gal(K))$ and $\rho_{\ell_2}(\Gal(K))$ do not have a common finite simple group as a quotient.
It is known that group theoretical independence implies independence (cf. \cite[Lemme 2]{bible}). The following
results have been established in the line of papers \cite{bible}, \cite{gp}, \cite{bgp}, \cite{cadoret}.

{\em Let $K$ be a field of characteristic $p\ge 0$ 
and $X/K$ a separated algebraic scheme. Let $\Ll'=\Ll\smallsetminus \{p\}$.
Let $q\in\Nn$ and let $\rho_\ell$ be the 
representation of $\Gal(K)$ on $H^q(X_{\ol{K}}, \Qq_\ell)$.
\begin{itemize}
\item[(I1)] If $K$ is a finitely generated field of characteristic zero, then there exists a finite extension
$E/K$ such that the restricted family $(\rho_\ell|\Gal(E))_{\ell\in\Ll}$ is independent.
\item[(I2)] If $K$ is a function field over an algebraically closed field, then there exists a finite separable extension
$E/K$ such that $(\rho_\ell|\Gal(E))_{\ell\in\Ll'}$ is group theoretically independent.
\end{itemize}}

(I1) in the important special case where $K$ is a number field was proved by Serre in \cite{bible}.
(I1) in the special case $\trdeg(K/\Qq)>0$ was proved by Gajda and the author 
in \cite{gp}, answering a question of Serre
(cf. \cite[Section 3.2]{bible}, \cite[10.1?]{serre1994})
and Illusie (cf. \cite[5.5]{illusie}). (I2) was proved by B\"ockle, Gajda and the author
in \cite{bgp} and independently by Cadoret and Tamagawa in \cite{cadoret}.

One can not replace the word ``independent'' by ``group
theoretically independent'' in (I1): 
If $K$ is a finitely generated field of characteristic zero, then $\dim_{\Qq_\ell}(H^2(\Pp^1_{\ol{K}}, \Qq_\ell))=1$ and the 
action of $\Gal(K)$ on $H^2(\Pp^1_{\ol{K}}, \Qq_\ell)$ is given by the inverse of the cyclotomic 
character $\varepsilon_{K, \ell}: \Gal(K)\to \Zz_\ell^\times$ over $K$, and it can easily be seen from the prime number theorem 
in arithmetic progressions that 
the family of cyclotomic 
characters $(\varepsilon_{K, \ell})_{\ell\in\Ll}$ is {\em not} group theoretically 
independent. Assume for the moment that $K$ is a finitely generated field without any restriction on its characteristic $p\ge 0$. 
In the light of these facts the following question comes up naturally.

 {\em Are the cyclotomic
characters the only obstruction to group theoretical independence, i.e. is there 
a finite separable 
extension $E/K$ such that $(\rho_\ell|\Gal(E(\mu_{\infty})))_{\ell\in\Ll'}$ is group theoretically
independent? }

If $p>0$, then the answer to this question is ``yes'', because then the prime field 
$\Ff$ of $K$ is finite and hence $\ol{\Ff}$ agrees with the field $\Ff(\mu_\infty)$ obtained from $\Ff$ by
adjoining all roots of unity; one can thus apply (I2) 
to $X_{\ol{\Ff}K}/\ol{\Ff} K$ to conclude that there is a finite separable extension $E/K$ such
that $(\rho_{\ell}|\Gal(\ol{\Ff} E)_{\ell\in \Ll})$ is 
group theoretically independent.
Somewhat surprisingly the answer to the above question is ``no'' in the case
$p=0$ (cf. Corollary \ref{cmcoro}, Remark \ref{cmrema}), as we shall see in Appendix \ref{app}. 
One can construct counterexamples from
certain CM abelian varieties. We do have the
following affirmative result, however.

\begin{mthm} \label{thmA} (cf. Corollary \ref{indepcoro1}) Let $K$ be a finitely generated
field of characteristic $0$ 
and $X/K$ a smooth projective variety. Let $q\in\Nn$ and let $\rho_\ell$ be the 
representation of $\Gal(K)$ on $H^q(X_{\ol{K}}, \Qq_\ell)$. Then
there is a finite Galois extension $k/\Qq$ and a finite Galois extension $F/K$ such that 
the family $(\rho_\ell|\Gal(k_{\ab} F))_{\ell\in\Ll}$ is group theoretically independent, if $k_{\ab}$ stands for the maximal abelian Galois 
extension of $k$. Furthermore
$(\rho_\ell|\Gal(k_{\ab}^\dagger F))_{\ell\in\Ll}$ is group theoretically independent, if $k_{\ab}^\dagger$ is the compositum
$k_{\ab}^\dagger=\prod_{\ell\in\Ll} k_{ab}^{(\ell)}$ where $k_{ab}^{(\ell)}/k$ is the maximal abelian extension of $k$ which is 
unramified outside $\ell$ and of order prime to $\ell$. 
\end{mthm}
 
 In certain special cases it can be shown that group theoretical independence is achieved already over the smaller extension $F_{\cyc}^\dagger
 =\prod_{\ell\in\Ll} F(\mu_\ell)$ where $F(\mu_\ell)$ is the field obtained from $F$ in adjoining an $\ell$-th root of unity. 
 For example, if in the situation of the main theorem $q=1$ and $X$ is an abelian variety with $\End_{\overline{K}}(X)=\Zz$ and of dimension $2$, $6$ or odd, then, based on Serre's open image theorem \cite[2.2.8]{serre-cours1984}, one can see easily that in this case one can choose  the finite Galois extension
 $F/K$ in such a way that 
$(\rho_\ell|\Gal(F_{\cyc}^\dagger))_{\ell\in\Ll}$ is group theoretically independent. On the other hand, if $q=1$, $K$ is a number field and 
$X/K$ is an absolutely simple abelian variety with complex multiplication over $K$, then $(\rho_\ell|\Gal(F_{\cyc}^\dagger ))_{\ell\in\Ll}$ is not group theoretically independent 
for every finite Galois extension $F/K$ (cf. Corollary \ref{cmcoro}, Remark \ref{cmrema}).

There are three main ingredients to the proof: (1) We make strong use of certain concepts from group theory (profinite and algebraic). In particular we crucially use 
information about the structure of subgroups of $GL_n(\Ff_\ell)$ (cf. \cite{larsenpink} and \cite[Section 3]{bgp}) and about point groups of reductive algebraic groups defined over finite fields (cf. 
Proposition \ref{gtprop}). (2) We certainly use information about the ramification of the representations under consideration provided by constructability and semistability theorems in \'etale cohomology (cf. \cite{illusie}). Furthermore we make use of Caruso's solution of Serre's tame inertia conjecture (cf. \cite{caruso}) 
in order to control the ramification of $\rho_\ell$ ``at primes above $\ell$''. (3) Finally we invoke finiteness results for \'etale fundamental groups from geometric class field theory
(cf. \cite{katzlang} and \cite[Section 2]{gp}).

\begin{center}{\bf Notation}\end{center} Throughout this manuscript $\Ll$ denotes the set of all prime numbers. 
If $K$ is a field, then we denote by $\ol{K}$ an algebraic closure of $K$ and by $\Gal(K)$ its absolute Galois group. 
We denote by $K_{\ab}$ the maximal abelian extension of $K$ and by $K_{\solv}$ the 
maximal prosolvable extension of $K$. Furthermore we denote by $K_{\cyc}=K(\mu_\infty)$ the extension obtained from
$K$ by adjoining all roots of unity and put $K_{\cyc}^\dagger=\prod_{\ell\in\Ll} K(\mu_\ell)$. If $K$ is a number field,
then for $\ell\in\Ll$ we denote by $K_{\ab}^{(\ell)}$ the maximal abelian extension of $K$ 
which is unramified outside $\ell$ and of degree prime to $\ell$. Furthermore we define 
$K_{\ab}^\dagger:=\prod_{\ell\in \Ll} K_{\ab}^{(\ell)}$.

A {\em $K$-variety} is a separated algebraic 
$K$-scheme which is reduced and irreducible. 
For a profinite group $G$ we shall denote by $S_\ell(G)$ the 
normal subgroup of $G$ generated by its $\ell$-Sylow subgroups. If $\ell$ is clear from the context, then we 
write $G^+$ instead of $S_\ell(G)$. 

\section{Preliminaries on group theory}\label{gtSec}
This section is devoted to the concepts from group theory used in this paper.
The following proposition about point groups of connected algebraic groups over finite fields is 
probably well-known to the expert. It will be applied later in the situation of the Main Theorem 
to suitable reductive envelopes of the images of the semisimplified mod-$\ell$ representations under 
consideration.

\begin{prop} \label{gtprop}
Let $F$ be a finite field of characteristic $\ell$ and $\ul{G}$ a connected 
algebraic group over $F$. 
Then the group $\ul{G}(F)/\ul{G}(F)^+$ is an {\em abelian} group of order prime to $\ell$.
\end{prop}

{\em Proof.} Let $\ul{U}$ be the unipotent radical of $\ul{G}$ and 
$\ul{H}=\ul{G}/\ul{U}$. Let $\ul{S}=[\ul{H}, \ul{H}]$ be the derived group of $\ul{H}$.  
The group $\ul{U}(F)$ is a normal subgroup of $\ul{G}(F)$ of $\ell$-power order. Hence 
$\ul{U}(F)\subset \ul{G}(F)^+\subset \ul{G}(F)$. Furthermore $\ul{H}(F)=\ul{G}(F)/\ul{U}(F)$ because
$\ul{U}$ is connected. It follows that $H(F)^+=\ul{G}(F)^+/\ul{U}(F)$ and $\ul{G}(F)/\ul{G}(F)^+
\cong \ul{H}(F)/\ul{H}(F)^+$. Thus we can assume that $\ul{G}$ is reductive (i.e. that $\ul{U}$ is trivial 
and $\ul{G}=\ul{H}$). Then
$\ul{T}:=\ul{G}/\ul{S}$ is a torus and $\ul{G}(F)/\ul{S}(F)$ embeds into $\ul{T}(F)$. Hence
$\ul{G}(F)/\ul{S}(F)$ is abelian of order prime to $\ell$. This implies that $\ul{G}(F)^+=\ul{S}(F)^+$. 

Let $\ul{Z}$ be the center of $\ul{G}$. Let $\ul{\widetilde{S}}\to \ul{S}$ be the simply connected 
covering of the semisimple group $\ul{S}$ and let $I$ be the image of $\ul{\widetilde{S}}(F)\to \ul{S}(F)$. Then
$\ul{S}(F)^+=I$ by a theorem of Steinberg (cf. \cite[12.4, 12.6]{steinberg1968}).
It thus suffices to show that $[\ul{G}(F), \ul{G}(F)]\subset I$. To prove this, let $a_1, a_2\in \ul{G}(F)$ and $s=[a_1, a_2]$. Then $s\in \ul{S}(F)$
and we denote by $\ol{a_i}$ the image of $a_i$ in $\ul{G}/\ul{Z}(F)$. We have central isogenies
$$\ul{\widetilde{S}}\buildrel f \over\longrightarrow \ul{S}\buildrel g \over\longrightarrow \ul{G}/\ul{Z}.$$
Choose for $i\in\{1, 2\}$ an element $\tilde{a_i}\in \ul{\widetilde{S}}(\ol{F})$ such that $g\circ f(\tilde{a_i})=\ol{a_i}$ and 
define $\tilde{s}=[\tilde{a}_1, \tilde{a}_2]$.
There exist elements $z_i\in \ul{Z}(\ol{F})$ such that $f(\tilde{a}_i)=z_i a_i$. Thus $$f(\tilde{s})=[f(\tilde{a}_1), 
f(\tilde{a}_2)]=[z_1 a_1, z_2 a_2]=[a_1, a_2]=s.$$
For every $\sigma\in \Gal(F)$ and every $i\in\{1, 2\}$ there exists an element $k_i\in \ker(g\circ f)(\ol{F})$ such that 
$\tilde{a_i}^\sigma=k_i \tilde{a_i}$. Hence 
$\tilde{s}^\sigma=[k_1 \tilde{a}_1, k_2 \tilde{a}_2]=\tilde{s},$
because the $k_i$ lie in the center of $\ul{\widetilde{S}}(\ol{F})$. It follows that 
$\tilde{s}\in \ul{\widetilde{S}}(F)$ is 
$F$-rational. Thus $s\in I$ as desired.\hfill $\Box$

\begin{defi} \label{classes} For $d\in\Nn$ we shall denote by $\mf{B}(d)$ the class of all finite groups of order $\le d$, and by 
 $\mf{Jor}(d)$ the class of all finite groups $G$ which admit an abelian normal subgroup $N$
such that $G/N\in \mf{B}(d)$. For $\ell\in\Ll$ denote by $\mf{Lie}_\ell$ the class of all finite simple groups of Lie type in characteristic
$\ell$.
Let $\mathfrak{Lie}_\ell(d)$ be the class of all finite groups
$G$ which admit an abelian normal subgroup $N$ such that $|N|$ is coprime to $\ell$, $|N|\le d$ and $G/N$ is isomorphic
to a finite product of groups in $\mf{Lie}_\ell$. 
\end{defi}

Note that the product is allowed to be empty and 
thus the trivial group lies in $\mf{Lie}_\ell(d)$ (but not in $\mf{Lie}_\ell$). 
The following theorem about finite subgroups of $\GL_n(\ol{\Ff}_\ell)$ is a corollary to a result of Larsen and Pink  (cf. \cite{larsenpink}) which was established in
\cite[Section 3]{bgp}. It slightly generalizes \cite[Thm. 3']{bible} and \cite[Thm. 4]{bible}. We shall apply it in order to understand the group theoretical properties 
of the images of certain mod-$\ell$ representations later.

\begin{thm} \label{nori} (cf. \cite[Section 3]{bgp}) For every $n\in\Nn$ there exists a constant $J'(n)\ge 5$ with the 
following property: For every $\ell\in\Ll$ and every finite
subgroup $G$ of $\GL_n(\ol{\Ff}_\ell)$ the group $G/G^+$ lies in $\mf{Jor}(J'(n))$.
Moreover, if $\ell> J'(n)$ and $P$ denotes the maximal normal $\ell$-subgroup of $G^+$, then
$G^+/P$ lies in $\mf{Lie}_\ell(2^{n-1})$. 
\end{thm}

Let $G$ be a profinite group and $H$ a finite simple (not necessarily non-abelian) group. We call $H$ a 
{\em Jordan-H\"older factor of $G$} if there exists a closed normal subgroup $G_1$ of $G$ and an open 
normal subgroup $G_2$ of $G_1$
and a continuous 
isomorphism $G_1/G_2\cong H$. We denote by $\JH(G)$ the class of all Jordan-H\"older factors of $G$. Let $\FSQ(G)$ be the class of all finite simple (not necessarily non-abelian) quotients of $G$. Then
$\FSQ(G)\subset \JH(G)$. The proof of the following elementary Lemma is left to the reader.

\begin{lemm}\label{JHexact}
\begin{enumerate}
\item[(a)] If $1\to G'\to G\buildrel \pi\over \longrightarrow 
G''\to 1$ is an exact sequence of profinite groups, then
$\JH(G)= \JH(G')\cup \JH(G'')$ and $\FSQ(G)\subset \FSQ(G')\cup \FSQ(G'')$.
\item[(b)] Let $G$ be a profinite group. Let $(N_i)_{i\in I}$ be a family of closed normal subgroups of $G$
and $N=\bigcap_{i\in I} N_i$.
Then $\JH(G/N)\subset \bigcup_{i\in I} \JH(G/N_i)$.
\end{enumerate}
\end{lemm}

Theorem \ref{nori} will frequently enter our considerations through the following remark.

\begin{rema} \label{JHrema}
{\em Let $G$ be a finite subgroup of $\GL_n(\Ff_\ell)$ where $\ell\ge J'(n)$. Let $P$ be the maximal normal $\ell$-subgroup of $G^+$.
\begin{enumerate}
\item[a)] The group $G/G^+$ lies in $\mf{Jor}(J'(n))$ (cf. Theorem \ref{nori}) and thus  $$\JH(G/G^+)\subset \mf{B}(J'(n))\cup \{\Zz/p: p\in \Ll\setminus \{\ell\}\}$$ by Lemma
\ref{JHexact}.
\item[b)] The group $G^+/P$ lies in $\mf{Lie}_\ell(2^{n-1})$ (cf. Theorem \ref{nori}) and thus $$\JH(G^+)\subset 
\mf{Lie}_\ell \cup \{\Zz/p: p\in \Ll\smallsetminus \{\ell\}, p\le 2^{n-1}\}\cup \{\Zz/\ell\}$$ by Lemma \ref{JHexact}. 
Furthermore $\FSQ(G^+)$ cannot contain groups of order prime to $\ell$ as $G^+$ is generated by its $\ell$-Sylow subgroups, 
and consequently $\FSQ(G^+) \subset \mf{Lie}_\ell\cup \{\Zz/\ell\}$. 
\end{enumerate}}
\end{rema}





For technical reasons the following Lemma will be useful.

\begin{lemm} \label{goingup} Let $G$ be a finite subgroup of $\GL_n(\Ff_\ell)$. Assume that $\ell>J'(n)$. Let $N$ be a normal subgroup of $G^+$. 
If $\JH(G/N)\cap \Lie_\ell=\emptyset$, then $N=N^+$. 
\end{lemm}

{\em Proof.} Let $P$ be the maximal normal $\ell$-subgroup of $G$ and consider the exact sequence of groups
$$1\to N/P\cap N\to G^+/P\to G^+/NP\to 1.$$
The group $G^+/NP$ is a quotient of a $G^+/N$ which is in turn isomorphic to a normal subgroup of $G/N$. 
Hence, by Lemma \ref{JHexact}, $\JH(G^+/NP)\cap\Lie_\ell=\emptyset$. On the other hand
$G^+/NP$ is a quotient of $G^+/P$. As $\ell>J'(n)$ we have 
$$\FSQ(G^+/NP)\subset \FSQ(G^+/P)\subset \Lie_\ell$$ (cf. Remark \ref{JHrema}).  
It follows that $\FSQ(G^+/NP)=\emptyset$, hence $G^+/NP$ is the trivial group. By the exact sequence above 
$N/P\cap N \cong  G^+/P$ is a group which is generated by its $\ell$-Sylow subgroups. As $P\cap N$ is an $\ell$-group, it follows
that $N$ is generated by its $\ell$-Sylow subgroups as well, i.e. $N=N^+$.\hfill $\Box$

\section{The monodromy groups of the mod-$\ell$ representations}\label{dfSec}
Let $K$ be a field of characteristic zero and $X/K$ a smooth projective variety. Fix $q\in \Nn$.
Let $\rho_\ell$ (resp. 
$\ol{\rho}_\ell$) be the representation of $\Gal(K)$ on $H^q(X_{\ol{K}}, \Qq_\ell)$ (resp. on
$H^q(X_{\ol{K}}, \Ff_\ell)$). The following Lemma explains the relation between the $\ell$-adic monodromy groups 
$\rho_\ell(\Gal(K))$ and the mod-$\ell$ monodromy groups $\orho_\ell(\Gal(K))$. The rest of this section is then 
devoted to the images of the mod-$\ell$ representations.

\begin{lemm} \label{gabber} 
There is a constant $D$ such that for every prime number $\ell\ge D$ there
is an epimorphism $\pi_\ell: \rho_\ell(\Gal(K))\to \ol{\rho}_\ell(\Gal(K))$ such that $\pi_\ell\circ \rho_\ell=\ol{\rho}_\ell$ and 
such that $P_\ell:=\ker(\pi_\ell)$ is a pro-$\ell$ group. 
\end{lemm}

{\em Proof.} By a theorem of Gabber (cf. \cite{gabber1983}) there exists
a constant $D$ such that $H^q(X_{\ol{K}}, \Zz_\ell)$ and $H^{q+1}(X_{\ol{K}}, \Zz_\ell)$ torsion
free $\Zz_\ell$-modules for all prime numbers $\ell\ge D$. 
For every $\ell\in\Ll$ there is an exact sequence of $\Gal(K)$-modules
$$\begin{array}{ll}
 & H^q(X_{\ol{K}}, \Zz_\ell)\buildrel \ell\over \longrightarrow H^q(X_{\ol{K}}, \Zz_\ell)\to H^q(X_{\ol{K}}, \Ff_\ell)\to \\
 \to & H^{q+1}(X_{\ol{K}}, \Zz_\ell)\buildrel \ell\over \longrightarrow H^{q+1}(X_{\ol{K}}, \Zz_\ell).\end{array}$$

We put $V_\ell:=H^q(X_{\ol{K}}, \Qq_\ell)$, $T_\ell:=H^q(X_{\ol{K}}, \Zz_\ell)$ and $W_\ell:=H^q(X_{\ol{K}}, \Ff_\ell)$.
Then, for all $\ell\in\Ll$ with $\ell>D$ the natural map 
$T_\ell \otimes \Ff_\ell\to W_\ell$
is an isomorphism, because $H^{q+1}(X_{\ol{K}}, \Zz_\ell)[\ell]=0$. Furthermore 
$V_\ell=T_\ell\otimes \Qq_\ell$ (by definition), and the natural map $T_\ell\to V_\ell$ must be injective because $T_\ell$ is
torsion free. 
We denote by $\rho_\ell'$ the representation
of $\Gal(K)$ on the finitely generated free $\Zz_\ell$-module $T_\ell$. The canonical maps
$$\GL_{W_\ell}(\Ff_\ell)\buildrel F_\ell \over \longleftarrow\GL_{T_\ell}(\Zz_\ell)\buildrel G_\ell \over
\longrightarrow
\GL_{V_\ell}(\Qq_\ell)$$
induce by restriction epimorphisms
$$\ol\rho_\ell(\Gal(K))\buildrel f_\ell \over\longleftarrow \rho_\ell'(\Gal(K))\buildrel g_\ell \over
\longrightarrow \rho_\ell(\Gal(K))$$
such that $f_\ell\circ \rho_\ell'=\ol{\rho}_\ell$ and $g_\ell\circ \rho_\ell'=\rho_\ell$. Furthermore
$\ker(f_\ell)$ is pro-$\ell$ because $\ker(F_\ell)$ is pro-$\ell$, and $g_\ell$ is injective because $G_\ell$ is injective.
It follows that $\pi_\ell:=f_\ell\circ (g_\ell)^{-1}$ is an epimorphism $\rho_\ell(\Gal(K))\to \ol\rho_\ell(\Gal(K))$ 
such that $\rho_\ell\circ \pi_\ell=\ol{\rho}_\ell$, and
$P_\ell:=\ker(\pi_\ell)$ is pro-$\ell$. 
\hfill $\Box$

For a non-archimedian place $v$ of a number field $K$ we shall denote by $\chara(v)$ its residue characteristic.
Recall the definition of $K_{\ab}^\dagger$ from the notation section.

\begin{lemm} \label{huilemm} 
Let $K$ be a number field, $X/K$ a smooth projective geometrically irreducible variety and $q\in\Nn$. Let
$\ol{\rho}_\ell$ be the representation of $\Gal(K)$ on $H^q(X_{\ol{K}}, \mathbb{F}_\ell)$.
Then there exists a finite Galois extension $E/K$ such that $\ol{\rho}_\ell(\Gal(E_{\ab}^\dagger))\subset \ol{\rho}_\ell(\Gal(K))^+$
for every $\ell\in\Ll$.
\end{lemm}

{\em Proof.} We denote by $\rho_\ell$ the representation of $\Gal(K)$ on $H^q(X_{\ol{K}}, \mathbb{Q}_\ell)$.
By \cite[Cor. 2.3]{illusie} there exists 
a finite set $S$ of places of $K$ such that for every non-archimedian place $v$ 
of $K$ outside $S$ and every prime number $\ell\neq \chara(v)$ the representation ${\rho}_\ell$ is unramified at $v$. 
Furthermore there exists a finite Galois extension of $K'/K$, such that for every non-archimedian place $v'$ of $K'$ above 
$S$ and every prime number $\ell\neq \chara(v')$ the group ${\rho}_\ell(I_{v'})$ is a pro-$\ell$ group (cf. \cite{dejong},
\cite[6.3.2]{berthelot}). After replacing $K'$ by a larger finite Galois extension of $K$ we can assume in addition that 
$\ol{\rho}_\ell(\Gal(K'))=\{e\}$ for all $\ell\le D$ where $D$ is the constant from Lemma \ref{gabber}. It follows (via Lemma \ref{gabber})
that for every place $v'$ of $K'$ and every $\ell\neq \chara(v')$
the group $\ol{\rho}_\ell(I_{v'})$ is an $\ell$-group, which is trivial if $v'$ does not lie over $S$.

Let $W_\ell^{ss}$ be the semisimplification of the $\Gal(K)$-module  $H^q(X_{\ol{K}}, \mathbb{F}_\ell)$ and $\ol{\rho}_\ell^{ss}$ be the representation of $\Gal(K)$ on $W_\ell^{ss}$. By the above we see that for
every place $v'$ of $K'$ and every prime number $\ell\neq \chara(v')$ the group
$\ol{\rho}_\ell^{ss}(I_{v'})$ is an $\ell$-group, which is trivial if $v'$ does not lie over $S$. Let $\ell_0=\max(\chara(w): w\in S)$. Let $\ell\in \Ll$ with $\ell>\ell_0$. If
$w'$ is a place of $K'$ with $\chara(w')=\ell$, then by Caruso \cite[Thm. 1.2]{caruso} the 
weight of the tame inertia group $I^t(w')$ acting via the 
contragredient of the semisimplification of the restricted representation $\orho_\ell^{ss}|I(w')$ are comprised in 
the inval $[0, eq]$ where $e$ is the ramification index of $K'_{w'}/\Qq_\ell$. Furthermore $\dim(W_\ell^{ss})$ does not depend on  $\ell$. 
Altogether we see that $\orho_\ell^{ss}|\Gal(K')$ satisfies the conditions (a)-(d) of \cite[Section 3.3]{wintenberger2002}. 
By \cite[Thm. 4]{wintenberger2002} (or by \cite[Thm. 2.3.5]{hui})
there exists, after replacing $\ell_0$ by a larger constant, a finite Galois
extension $L/K'$ and for every prime number $\ell\ge \ell_0$ 
a reductive algebraic subgroup $\ul{G}_\ell/\Ff_\ell$ of 
$\ul{GL}_{W_\ell^{ss}}$ with the following properties:
\begin{enumerate}
\item[(1)] $\ol{\rho}_\ell^{ss}(\Gal(L)) 
 \subset \ul{G}_\ell(\Ff_\ell)$ 
for every prime number $\ell\ge \ell_0$. 
\item[(2)] $\ol{\rho}_\ell^{ss}(\Gal(L)) ^+=\ul{G}_\ell(\Ff_\ell)^+$ for every prime number $\ell\ge \ell_0$. 
\end{enumerate}
We can replace $L$ by its Galois closure and enlarge $\ell_0$ accordingly in order to assume $L/K$ is Galois.
It follows from (1), (2) and Proposition \ref{gtprop} that $\ol{\rho}_\ell^{ss}(\Gal(L))/
\ol{\rho}_\ell^{ss}(\Gal(L))^+$ is abelian for every prime number 
$\ell\ge \ell_0$. The kernel $P_\ell$ of the natural epimorphism $g: \ol{\rho}_\ell(\Gal(L))\to 
\ol{\rho}_\ell^{ss}(\Gal(L))$ is an $\ell$-group; hence it lies in $\ol{\rho}_\ell(\Gal(L))^+$. Thus
$g$ induces an isomorphism 
$$\ol{\rho}_\ell(\Gal(L))/\ol{\rho}_\ell(\Gal(L))^+\cong 
\ol{\rho}_\ell^{ss}(\Gal(L))/\ol{\rho}_\ell^{ss}(\Gal(L))^+.$$
It follows that $\ol{\rho}_\ell(\Gal(L))/\ol{\rho}_\ell(\Gal(L))^+$ is abelian for every $\ell \ge \ell_0$.
In particular
 $$\ol{\rho}_\ell(\Gal(E_{\ab}))\subset \ol{\rho}_\ell(\Gal(L))^+\subset \ol{\rho}_\ell(\Gal(K))^+$$
 for every finite extension $E/L$ and every $\ell\ge \ell_0$.
We now choose $E$ to be a finite Galois extension of $K$ containing $L\cdot \prod_{\ell\le \ell_0} K(\ol{\rho}_\ell)$. 
Then $\ol{\rho}_\ell(\Gal(E_{\ab}))\subset  \ol{\rho}_\ell(\Gal(K))^+$ for every $\ell\in \Ll$. Moreover, for every $\ell\in \Ll$, the 
group
$\ol{\rho}_\ell(I_v)$ is an $\ell$-group and hence contained in $\ol{\rho}_\ell(\Gal(K))^+$ for every place $v$ of $E$ with $\chara(v)\neq \ell$.
Thus $\ol{\rho}(\Gal(E_{ab}^{(\ell)}))\subset \ol{\rho}_\ell(\Gal(K))^+$ for all $\ell\in\Ll$. 
Thus the assertion follows with the above choice of $E$.\hfill $\Box$

We shall now generalize Lemma \ref{huilemm} to the situation where the ground field $K$ is
an arbitrary finitely generated extension of $\Qq$, possibly of transcendence degree $\ge 1$. For this we use
a specialization argument along with a finiteness theorem for unramified Jordan extensions from \cite{gp}; this theorem from \cite{gp}  in turn
relies on a finiteness theorem in geometric class field theory of Katz and Lang (cf. \cite{katzlang}) and on some finiteness results for 
geometric fundamental groups from SGA.

\begin{prop} \label{divisionfields}
Let $K/\Qq$ be a finitely generated extension of fields. Let $X/K$ be a smooth projective variety. For 
every prime number $\ell$ let $\ol{\rho}_\ell$ 
be the representation of $\Gal(K)$ on $H^q(X_{\ol{K}}, \Ff_\ell)$. 
There exists a finite Galois extension $E/K$ and a finite Galois extension $k/\Qq$ such that
$\ol{\rho}_\ell(\Gal(k^\dagger_{\ab} E))\subset \ol{\rho}_\ell(\Gal(K))^+ $ for every $\ell\in\Ll$. 
\end{prop}

{\em Proof.} There exists a finite Galois extension $K'/K$ such that $X_{K'}$ splits up into a disjoint 
sum of geometrically connected (smooth projective) $K'$-varieties. Once the proposition is true for every connected component of 
$X_{K'}/K'$ it will follow for $X/K$. We may thus assume right from the outset that $X/K$ is geometrically connected.

There exists  a $\Qq$-variety $S$ with function field $K$. Moreover, by the usual spreading-out principles, there
exists after replacing $S$ by one of its dense open subschemes a smooth projective morphism  
$f: \mathscr{X}\to S$ with generic fibre $X$.
The stalk of $f_* \mathscr{O}_{\mathscr{X}}$ at the generic point of $S$ is zero because $X/K$ is 
geometrically connected. Now, after replacing $S$ by one of its non-empty open subschemes and 
shrinking $\mathscr{X}$ accordingly, we may assume that $f_* \mathscr{O}_{\mathscr{X}}=0$, where $\mathscr{O}_{\mathscr{X}}$ stands for
the structure sheaf of $\mathscr{X}$. Then $f$ has geometrically connected fibres (c.f. \cite[4.3.4]{EGAIII}).
Shrinking $S$ and $\mathscr{X}$ once more one can assume that the \'etale sheaves $R^q f_* \Ff_\ell$
are lisse and compatible with any base change \cite[Cor. 2.2]{illusie}. In particular the representation $\ol{\rho}_\ell$ factors 
through\footnote{We take the \'etale fundamental group with resprect to the geometric generic base point of $S$ afforded by the choice of $\ol{K}$.}  the
\'etale fundamental group $\pi_1(S)$. Let $G_\ell=\ol{\rho}_\ell(\pi_1(S))$ and let $n$ be an upper bound for the dimensions of the $\Ff_\ell$-vector spaces $H^q(X_{\ol{K}}, \Ff_\ell)$ 
($\ell\in\Ll$). The existence of such an upper bound is guaranteed by \cite[Thm. 1.1]{illusie}. Now $\orho_\ell(\Gal(K))$ is isomorphic to a subgroup of $\GL_n(\Ff_\ell)$ for
all $\ell\in\Ll$. Thus there exists a
constant $J'(n)$ such that $G_\ell/G_\ell^+\in \mf{Jor}(J'(n))$ for all $\ell\in\Ll$ (cf. Theorem \ref{nori}). Hence, by \cite[Prop. 2.2]{gp}, there exists an open
normal subgroup $U$ of $\pi_1(S)$ such that $\ol{\rho}_\ell(U\cap \pi_1(S_{\ol{\Qq}}))\subset G_\ell^+$ for all $\ell\in\Ll$. Let $s\in S$ be a closed point.
Let $S'$ be the finite \'etale Galois cover of $S$ corresponding to $U$  and pick a closed point $s'\in S'$ over $s$. Let $E$ be the function field of $S'$.

Let us consider the following commutative diagram of profinite groups:
$$\xymatrix{ & & G_\ell & \\
1 \ar[r]& \pi_1(S_{\ol{\Qq}})\ar[r] & \pi_1(S)\ar[u]_{\ol{\rho}_\ell}\ar[r] & \Gal(\Qq)   \\
 & & &\Gal(k(s))\ar@{->}[u]\ar[ul]_{s_*} \\
   }$$
(The map $s_*$ is well-defined only up to conjugation.) If $X_s=\mathscr{X}\times_S\Spec (k(s))$ is the special fibre of $\mathscr{X}$ over $S$, then, by the base change compatibility alluded to above, the representation 
$\ol{\rho}_\ell \circ s_*$ of $\Gal(k(s))$ on $H^q(X_{\ol{K}}, \Ff_\ell)$ is isomorphic to the representation of $\Gal(k(s))$ on
$H^q(X_{s, \ol{k(s)}}, \Ff_\ell)$. Furthermore $X_s$ is a smooth projective geometrically connected variety over the number field $k(s)$. By Lemma
\ref{huilemm} there is a finite Galois extension $k/\Qq$ containing $k(s)$ such that $$\ol{\rho}_\ell\circ s_*(\Gal(k_{\ab}^\dagger))\subset \ol{\rho}_\ell\circ s_*(\Gal(k))^+\subset 
G_\ell^+.$$
After replacing $k$ by a finite extension we can also assume that $k\supset k(s')$. 
Now there is a commutative diagram with exact rows

$$\xymatrix{ & & G_\ell & \\
1 \ar[r]& \pi_1(S_{\ol{\Qq}})\ar[r] & \pi_1(S)\ar[u]_{\ol{\rho}_\ell}\ar[r] & \Gal(\Qq)   \\
 & & &\Gal(k(s))\ar@{->}[u]\ar[ul]_{s_*} \\
 1 \ar[r]& \pi_1(S'_{\ol{\Qq}})\ar[r]\ar@{->}[uu] & \pi_1(S'_{k_{\ab}^\dagger})\,\ar@{->}[uu]\ar@<2pt>[r] & \Gal(k_{\ab}^\dagger)\ar@{->}[u]\ar[r]\ar@<2pt>[l] &1.}$$

We already know that $\ol{\rho}_\ell(\pi_1(S'_{\ol{\Qq}}))$ and $\ol{\rho}_\ell(s_*(\Gal(k_{\ab}^\dagger)))$ are contained in $G_\ell^+$. As $\pi_1(S'_{k_{\ab}^\dagger})$ is generated
by $\pi_1(S'_{\ol{\Qq}})$ and $s_*(\Gal(k_{\ab}^\dagger))$ we conclude that $\ol{\rho}_\ell(\Gal({k_{\ab}^\dagger E}))=\ol{\rho}_\ell(\pi_1(S'_{k_{\ab}^\dagger}))\subset G_\ell^+$
as desired. \hfill $\Box$

\begin{rema} In the situation of Proposition \ref{divisionfields} it is easy to see (with the help of Remark \ref{JHrema} and Lemma \ref{JHexact}) that 
$$\JH(\orho_\ell(\Gal(k_{\ab}^\dagger E)))
\subset \Lie_\ell\cup \{\Zz/p: p\in\Ll\setminus \{\ell\}, p\le 2^{n-1}\}\cup \{\Zz/\ell\}$$
for all but finitely many primes $\ell\in\Ll$. 
\end{rema}

Note, however, that this does not rule out the possibility that for some small prime number $p$  there exist infinitely many
$\ell\in\Ll$ such that 
$\Zz/p$ is a finite simple quotient of $\orho_\ell(\Gal(k_{\ab}^\dagger E))$. Hence Proposition \ref{divisionfields} alone does not imply group theoretical independence for the family
 $(\orho_\ell|\Gal(k_{\ab}^\dagger E))_{\ell\in\Ll}$; we need additional arguments to establish the Main Theorem.

\section{Independence results}\label{ISec}

In the following theorem we shall prove among other things that in the situation of Proposition \ref{divisionfields} one can, after replacing 
$E$ by a finite extension $F$ which is Galois over $K$, achieve a very good control over the possible finite simple quotients of 
$\orho_\ell(\Gal(k_{\ab}^\dagger F))$ and of $\rho_\ell(\Gal(k_{\ab}^\dagger F))$. In particular we shall see that for a suitable choice of $F$
the groups $\orho_\ell(\Gal(k_{\ab}^\dagger F))$ cannot have a finite simple quotient of order prime to $\ell$ any more.
The argument is of a group
theoretical nature. The Main Theorem about group theoretical independence along with
some variants will then follow quite easily. 

\begin{thm} \label{images} Let $K/\Qq$ be a finitely generated field extension. 
Let $X/K$ be a smooth projective
variety. Fix $q\in\Nn$. Let $\rho_\ell$ (resp. 
$\ol{\rho}_\ell$) be the representation of $\Gal(K)$ on $H^q(X_{\ol{K}}, \Qq_\ell)$ (resp. on
$H^q(X_{\ol{K}}, \Ff_\ell)$). Let $\ell_0\in\Nn$. 
Then there is a finite Galois extension $F/K$ and 
a finite Galois extension $k/\Qq$ with the following properties.
\begin{enumerate}
\item[(a)] For every algebraic extension $\Omega/F$ and every $\ell\le \ell_0$ in $\Ll$ the group $\ol{\rho}_\ell(\Gal(\Omega))$ is trivial and the group
$\rho_\ell(\Gal(\Omega))$ is a pro-$\ell$ group.
\item[(b)] For every solvable Galois extension $\Omega/k_{\ab}^\dagger F$ and every $\ell\in\Ll$ we have
$\orho_\ell(\Gal(\Omega))=\orho_\ell(\Gal(\Omega))^+$, and there exists a closed normal pro-$\ell$ subgroup $Q_\ell$ of $\rho_\ell(\Gal(\Omega))$ such that 
$\rho_\ell(\Gal(\Omega))/Q_\ell\cong \orho_\ell(\Gal(\Omega))$. 
\item[(c)] For every  solvable Galois extension $\Omega/k_{\ab}^\dagger F$ and every $\ell\in\Ll$ we have $$\FSQ(\orho_\ell(\Gal(\Omega)))\subset \Lie_\ell\cup \Zz/\ell
\ \mbox{and}\ \FSQ(\rho_\ell(\Gal(\Omega)))\subset \Lie_\ell\cup \{ \Zz/\ell\}.$$ 
\end{enumerate}
\end{thm}

{\em Proof.} There exists $n\in\Nn$ such that $$\dim_{\Qq_\ell}(H^q(X_{\ol{K}}, \Qq_\ell)\le \dim_{\Ff_\ell}(H^q(X_{\ol{K}}, \Ff_\ell)\le n$$ for all $\ell\in\Ll$ (cf. \cite[Thm. 1.1]{illusie}).
We can assume right from the outset that $\ell_0\ge \max(J'(n), D)$ where $J'(n)$ is the constant from Theorem \ref{nori} and $D$ is the constant from
Lemma \ref{gabber}. We put $G_\ell:=\rho_\ell(\Gal(K))$ and $\ol{G}_\ell=\orho_\ell(\Gal(K))$. 
For every $\ell\in \Ll$ the maximal normal pro-$\ell$ subgroup $P_\ell$ of $G_\ell$ is open. 

By Proposition \ref{divisionfields} there exists a finite Galois extension $E/K$ such that $\orho_\ell(\Gal(k_{\ab}^\dagger E))\subset G_\ell^+$ for every $\ell\in\Ll$. 
Let $E':=E\cdot \prod_{\ell\le \ell_0}\ol{K}^{\rho_\ell^{-1}(P_\ell)}\cdot \ol{K}^{\ker(\orho_\ell)}$, pick a prime number $\ell_1>\max(\ell_0, [E':K])$ and let 
$F:=E'\cdot \prod_{\ell_0< \ell\le \ell_1} \ol{K}^{\ker(\orho_\ell)}$. The extensions $E'/K$ and $F/K$ are 
finite Galois extension because $P_\ell$ is open and normal in $G_\ell$ and $\ol{G}_\ell$ is finite.
Note that for every algebraic extension $\Omega/F$ the group $\rho_\ell(\Gal(\Omega))$ is pro-$\ell$ for every $\ell\le \ell_0$ in $\Ll$ and 
$\orho_\ell(\Gal(\Omega))$ is trivial for every $\ell\le \ell_1$ in $\Ll$. In particular a) holds true. 

To prove b) let $\Omega$ be a solvable Galois extension of 
$k_{\ab}^{\dagger} F$. We already know that b) holds true for all $\ell\le \ell_0$. For $\ell >\ell_0$, by Lemma \ref{gabber}, the group $\rho_\ell(\Gal(\Omega))$ is an extension
of $\orho_\ell(\Gal(\Omega))$ by a pro-$\ell$ group. As we know that $\orho_\ell(\Gal(\Omega))$ is trivial for all $\ell\le \ell_1$ it follows that $\rho_\ell(\Gal(\Omega))$ is pro-$\ell$ for
all $\ell\le \ell_1$. Thus b) holds true for every $\ell\le \ell_1$, and to establish b) completely it suffices to prove the following

{\bf Claim:} {\em $\orho_\ell(\Gal(\Omega))=\orho_\ell(\Gal(\Omega))^+$ for all $\ell>\ell_1$. }

Let $\mf{C}$ be the class of all prime cyclic groups. 
We shall now compute $\JH(\Gal(k_{\ab}^\dagger F/K))$ and then establish the claim with the help of Lemma \ref{goingup}.
By Lemma \ref{JHexact}
$$\JH(\Gal(k_{\ab}^\dagger F/K))\subset \JH(\Gal(k_{\ab}^\dagger K/K))\cup \JH(\Gal(E'/K))\cup \bigcup_{\ell_0<\ell\le \ell_1} \JH(\ol{G}_\ell),$$
and $\JH(\Gal(E'/K))\subset \mf{B}([E'/K])\subset \mf{B}(\ell_1)$. Moreover $\JH(\Gal(k_{\ab}^\dagger K/K))\subset \mf{C}$. If $\ell_0<\ell\le \ell_1$,
then $J'(n)<\ell$, and hence 
$$ \JH(\ol{G}_\ell)\subset \JH(\ol{G}_\ell^+)\cup \JH(\ol{G}_\ell/\ol{G}_\ell^+).$$
We have $\JH(\ol{G}_\ell^+)\subset \Lie_\ell\cup \mf{C}$ (cf. Remark \ref{JHrema}), and
$$\JH(\ol{G}_\ell/\ol{G}_\ell^+)\subset \mf{B}(J'(n))\cup \mf{C}\subset \mf{B}(\ell_1)\cup \mf{C}$$
(cf. Remark \ref{JHrema}).
Altogether it follows that 
$$\JH(\Gal(k_{\ab}^\dagger F/K))\subset \mf{C}\cup \mf{B}(\ell_1)\cup \bigcup_{\ell_0<\ell\le \ell_1} \mf{Lie}_\ell.$$
For $\ell>\ell_1$ the groups in $\mf{Lie_\ell}$ are non-commutative and generated by their $\ell$-Sylow subgroups, and hence they are neither contained in
$\mf{C}$ nor in  $\mf{B}(\ell_1)$. Together with a theorem of E. Artin about the orders of the finite simple groups of Lie type (cf. \cite[Thm. 5]{bible}, see also
\cite{artin}, \cite{KLMS}) we see that 
$$\JH(\Gal(k_{\ab}^\dagger F/K))\cap \mf{Lie}_\ell=\emptyset\quad \mbox{for all $\ell>\ell_1$}.$$

We now prove the claim. Let $\ell>\ell_1$. Let $N_\ell=\orho_\ell(\Gal(k_{\ab}^\dagger F))$ and $M_\ell=\orho_\ell(\Gal(\Omega))$. As $F\supset E$ by construction, we see that 
$N_\ell\subset G_\ell^+$. Moreover $G_\ell/N_\ell$ is a quotient of $\Gal(k_{\ab}^\dagger F/K)$. Hence $\JH(G_\ell/N_\ell)\cap \mf{Lie}_\ell=\emptyset$ for
all $\ell>\ell_1$. Lemma \ref{goingup} implies $N_\ell=N_\ell^+$. But then $M_\ell\subset N_\ell=N_\ell^+$, and moreover $N_\ell/M_\ell$ is a quotient 
of $\Gal(\Omega/k_{\ab}^\dagger F)$. As $\Omega/k_{\ab}^\dagger F$ is solvable we see that $\JH(N_\ell/M_\ell)\cap \mf{Lie_\ell}=\emptyset$ for all $\ell>\ell_1$. 
Applying Lemma \ref{goingup} once more, we see that $M_\ell=M_\ell^+$ for all $\ell>\ell_1$. This finishes up the proof of the claim and of Part (b). 

We now prove Part (c). If $\ell> \ell_1$, then (recall that $\ell_1\ge \ell_0\ge J'(n)$) $\FSQ(M_\ell) \subset \Lie_\ell\cup \{\Zz/\ell\}$ by Remark \ref{JHrema}.
If $\ell\le \ell_1$, then we even have $\FSQ(M_\ell)=\emptyset$ by part (a). Thus $\FSQ(M_\ell) \subset \Lie_\ell\cup \{\Zz/\ell\}$ for all $\ell\in\Ll$. As $\rho_\ell(\Gal(\Omega))$
is an extension of $M_\ell$ by a pro-$\ell$ group, the statement about $\rho_\ell$ in part (c) is also true. \hfill $\Box$ 

\begin{coro} \label{indepcoro1} Let $K/\Qq$ be a finitely generated field extension. 
Let $X/K$ be a smooth projective
variety. Fix $q\in\Nn$. Let $\rho_\ell$ (resp. 
$\ol{\rho}_\ell$) be the representation of $\Gal(K)$ on $H^q(X_{\ol{K}}, \Qq_\ell)$ (resp. on
$H^q(X_{\ol{K}}, \Ff_\ell)$).
Then there is a finite Galois extension $F/K$ and a finite Galois extension $k/\Qq$ such that 
for every solvable Galois extension $\Omega/k_{\ab}^\dagger F$ 
the families $(\rho_\ell|\Gal(\Omega))_{\ell\in\Ll}$  and $(\orho_\ell|\Gal(\Omega))_{\ell\in\Ll}$ are
group theoretically independent.
\end{coro}

{\em Proof.} Let $\ell_0=5$. Let $F/K$ and $k/\Qq$ be finite
Galois extensions such that the assertions (a) - (c) from Theorem \ref{images} hold true. From part (c) we have
$$\FSQ(\orho_\ell(\Gal(\Omega)))\subset \Lie_\ell\cup \{\Zz/\ell\}
\ \mbox{and}\ \FSQ(\rho_\ell(\Gal(\Omega)))\subset \Lie_\ell\cup \{\Zz/\ell\}$$
for all $\ell\in\Ll$. If $\ell\in \{2, 3\}$, then even 
$$\FSQ(\orho_\ell(\Gal(\Omega)))=\emptyset
\ \mbox{and}\ \FSQ(\rho_\ell(\Gal(\Omega)))\subset \{\Zz/\ell\}$$
by part (a). The class $\Lie_\ell$ consists of non-commutative groups, and $\Lie_{\ell_1}\cap \Lie_{\ell_2}=\emptyset$ for all 
prime numbers $5\le \ell_1<\ell_2$ by E. Artin's theorem (cf. \cite[Thm. 5]{bible}, see also
\cite{artin}, \cite{KLMS}). Hence $\FSQ(\orho_{\ell_1}(\Gal(\Omega)))\cap \FSQ(\orho_{\ell_2}(\Gal(\Omega)))=\emptyset$ for
all prime numbers $\ell_1\neq \ell_2$, and thus $(\orho_\ell)_{\ell \in \Ll}$ is group theoretically independent. Similarly the assertion about
$(\rho_\ell)_{\ell \in \Ll}$ is true. \hfill $\Box$

\begin{rema}\label{indepcoro2} If in the situation of Corollary \ref{indepcoro1} 
$$\Omega\in \{k_{\ab}^\dagger F, k_{\ab}F, k_{\solv} F, (kF)_{\ab}, (kF)_{\solv}\},$$
then $\Omega/k_{\ab}^\dagger F$ is solvable and hence
the families $$(\rho_\ell|\Gal(\Omega))_{\ell\in\Ll}\ \mbox{and}\ (\orho_\ell|\Gal(\Omega))_{\ell\in\Ll}$$ are
group theoretically independent.
\end{rema}

\begin{appendices}

\section{Appendix: A counterexample}\label{app}
The aim of this appendix is to prove that one cannot replace $k_{\ab}^\dagger F$ by $k_{\cyc}^\dagger F$ (or by 
$k_{\cyc} F$) in
Corollary \ref{indepcoro1}. In fact certain abelian varieties with complex multiplication over a number field
$K$ provide a counterexample. Throughout this section let $K$ be a number field and $A/K$ an absolutely simple
abelian variety with complex multiplication over $K$. We denote for every $m\in \Nn$ by $A[m]=\{x\in A(\overline{K}): mx=0\}$ the group of
$m$-torsion points of $A$ and for $\ell\in\Ll$ by $T_\ell(A)=\plim_{j} A[\ell^j]$ the $\ell$-adic Tate module of $A$.
Let $\eta_\ell$ (resp. $\ol{\eta}_\ell$) be the representation
of $\Gal(K)$ on $T_\ell(A)$ (resp. on $A[\ell]$). It is known that $\eta_\ell(\Gal(K))$ is abelian for all $\ell\in\Ll$.

\begin{thm}\label{cmdens} For every prime number $q$ there exists a set $S\subset \Ll$ of positive Dirichlet density such that
$\Zz/q\in \FSQ(\ol{\eta}_\ell(\Gal(K(\mu_\infty)))$ for every $\ell\in S$. 
\end{thm}

{\em Proof.} There exists a finite extension $K'/K$ such that $A_{K'}$ has good reduction everywhere
(cf. \cite[Thm. 7]{serretate1968}).

For every number field $F$ we let $\OO_F$ be its ring of integers and $\mathrm{Spl}(F)$ the set of
prime numbers $\ell$ which split completely in $F$. We always denote
by $\ul{T}_F=Res_{F/\Qq} \Gg_m$ the torus 
over $\Qq$ obtained as the Weil restriction of the multiplicative group over $F$. Furthermore we denote for every prime
$\mf{p}$ of $F$ by $F_{\mf p}$ the corresponding local field and by $U_{\mf p}(F)$ the group of units in the 
integer ring of $F_{\mf p}$. For any torus $\ul{T}/\Qq$, $\ell\in \Ll$ and $n\ge 0$ we define following
Ribet \cite[p. 77]{ribet1980}
$$\ul{T}(1+\ell^n \Zz_\ell):=\{x\in \ul{T}(\Qq_\ell): v_\ell(\chi(x)-1)\ge n\mbox{ for all } \chi\in\Hom_{\ol{\Qq}_\ell}(
\ul{T}_{\ol{\Qq}_\ell}, \Gg_{m, \ol{\Qq}_\ell})\}$$
where $v_\ell$ denotes the unique extension to $\ol{\Qq}_\ell$ of the canonical discrete valuation
of the complete field $\Qq_\ell$. We put $$\ul{T}(\Zz_\ell)=
\{x\in \ul{T}(\Qq_\ell): v_\ell(\chi(x)-1)\ge 0\mbox{ for all } \chi\in\Hom_{\ol{\Qq}_\ell}(
\ul{T}_{\ol{\Qq}_\ell}, \Gg_{m, \ol{\Qq}_\ell})\}$$
and $\ul{T}(\Ff_\ell)=
\ul{T}(\Zz_\ell)/\ul{T}(1+\ell \Zz_\ell)$. 
Note that $\ul{T}_F(\Qq_\ell)=\prod_{\mf p\mid \ell} F_{\mf p}^\times$ and $\ul{T}_F(\Zz_\ell)=
\prod_{\mf p\mid \ell} U_{\mf p}(F)$ (cf. \cite[Example 2.1]{ribet1980}). 
We denote by $I_F$ the idele group of $F$, define
$$I_F^{1}=\{(x_{\mf p})\in \prod_{\mf p} F_{\mf p}^\times: 
x_{\mf p}\in U_{\mf p}(F) \mbox{ for all } {\mf p} \nmid \infty\}$$
and view $\ul{T}_F(\Zz_\ell)$ as a subgroup of $I_F^1$. 

Let $I\subset \Ll$ be a finite set of prime numbers that contains the primes dividing 
$[\OO_E:E\cap \End_{K'}(A)]$ and
the primes that are ramified in $K'E$. Let $\ell\in\Ll\smallsetminus I$. Then
the $\OO_E\otimes \Zz_\ell$-module
$T_\ell(A)$ is free of rank $1$ (cf. \cite[Thm. 5]{serretate1968}) and hence
$$\eta_\ell(\Gal(K'))\subset \Aut_{\OO_E\otimes \Zz_\ell}(T_\ell(A))=(\OO_E\otimes \Zz_\ell)^\times=\ul{T}_E(\Zz_\ell).$$
Thus $\eta_\ell$ factors to a map $\Gal(K'_{\ab}/K')\to \ul{T}_E(\Zz_\ell)$ which is again denoted $\eta_\ell$. We compose $\eta_\ell$ with the Artin symbol $(-, K'_{\ab}/K')$
in order to obtain a map
$$\hat{\eta}_\ell: I_{K'}\to \Gal(K'_{\ab}/K')\buildrel{\eta_\ell}\over\longrightarrow \ul{T}_E(\Zz_\ell).$$ 
The image $U$ of $I_{K'}^1$ under the norm residue symbol is open in $\Gal(K_{\ab}'/K')$ and the fixed field of $U$ is the 
Hilbert class field $H$ of $K'$. If $\mf{p}$ is a prime of $K'$ and ${\mf p} \nmid \ell$, then $\eta_\ell$
is unramified at $\mf{p}$ by \cite[Thm 1]{serretate1968} because $A$ has good reduction everywhere. Hence $\hat{\eta}_\ell(U_\mf p)=\{e\}$. Thus $\hat{\eta}_\ell$ induces a map $\hat{\eta}'_\ell:
\ul{T}_{K'}(\Zz_\ell)=\prod_{\mf{p}\mid \ell} U_{\mf p} (K')\to \ul{T}_E(\Zz_\ell)$
and $\im(\hat{\eta}'_\ell)=\eta_\ell(\Gal(H))$. By a reformulation due to Serre and Tate (cf. \cite[Thm. 11, Cor. 2]{serretate1968}) of a theorem of Shimura and Taniyama there is a homomorphism $\psi: \ul{T}_{K'}\to \ul{T}_E$ of tori over $\Qq$ such that 
$\hat{\eta}'(x)=\psi_\ell(x^{-1})$ for all $x\in \ul{T}_{K'}(\Zz_\ell)$, where $\psi_\ell: \ul{T}_{K'}(\Zz_\ell)\to
\ul{T}_E(\Zz_\ell)$ is the homomorphism induced by $\psi$. Thus $\eta_\ell(\Gal(H))=\im(\psi_\ell)$. Let
$\ul{T}/\Qq$ be the image of $\psi_\ell$. Ribet proved that $\dim(T)\ge 2$ (cf. \cite[p. 87]{ribet1980}). 
(It is known that $\ul{T}$ agrees with the Mumford-Tate group of $A$, but
we will not need this fact.) Furthermore 
$\ol{\eta_\ell}(\Gal(H))=\im(\ol{\psi}_\ell)\subset \ul{T}(\Ff_\ell)$ if $\ol{\psi}_\ell: \ul{T}_{K'}(\Ff_\ell)\to \ul{T}_E(\Ff_\ell)$ is the homomorphism induced by $\psi$. 

By a theorem of Ribet (cf. \cite[Thm. 2.4]{ribet1980}) there is a constant $C$ such that $c_\ell:=
[\ul{T}(\Ff_\ell):
\ol{\eta}_\ell(\Gal(H))]\le C$ for all $\ell\in\Ll\smallsetminus I$. Furthermore, by \cite{bible},
there exists a finite extension $H'/H$ such that $(H'(\mu_{\ell^\infty}, A[\ell]))_{\ell\in\Ll}$ is a linearly disjoint sequence of 
extensions of $H'$. 
Now let $a\in \Nn$ be an exponent such that
$q^a>C[H':H]$. Let $L/\Qq$ be a Galois extension such that $\ul{T}\times\Spec(L)$ is a split torus and such that
$K'E(\mu_{q^a})\subset L$. Then $\mathrm{Spl}(L)\cap I=\emptyset$. Let $\ell\in 
\mathrm{Spl}(L)$. Then $\ul{T}\times {\Qq_\ell}$ is a split torus over $\Qq_\ell$
because $L$ can be embedded into $\Qq_\ell$.
Hence, if we put $d=\dim(\ul{T})$, then $\ul{T}(\Ff_\ell)=(\Ff_\ell^\times)^d$. If $\ell\in\mathrm{Spl}(L)$ does not
divide the polarization degree $\pi$ of $A$, then $H'(\mu_\ell)\subset H'(A[\ell])$ and 
$$[H'(\mu_\ell):H'][H'(A[\ell]):H'(\mu_\ell)] c_\ell=(\ell-1)^d,$$
$[H'(\mu_\ell):H']$ divides $\ell-1$, $c_\ell\le C[H':H]$, $q^a> C[H':H]$, $d\ge 2$ and 
$\ell=1\ \mathrm{mod}\ q^a$ because $\ell$ splits completely in $\Qq(\mu_{p^a})$. This forces
$[H'(A[\ell]):H'(\mu_\ell)]$ to be divisible by $q$. Furthermore $\Gal(
H'(\mu_{\ell^\infty})/H'(\mu_\ell))\cong \Zz_\ell$. It follows that $[H'(A[\ell], \mu_{\ell^\infty}): 
H'(\mu_{\ell^\infty})]$ is still divisible by $q$. Finally, by the linear disjointness of the sequence
$(H'(\mu_{\ell^\infty}, A[\ell]))_{\ell\in\Ll}$, one sees that $\prod_{\ell'\neq \ell} H'(\mu_{\ell'^\infty})$
is linearly disjoint from $H'(\mu_{\ell^\infty}, A[\ell])$ over $H'$ for every $\ell\in\Ll$. If follows that 
$|\ol{\eta}_\ell(\Gal(H'(\mu_\infty))|=[H'(\mu_\infty, A[\ell]):H'(\mu_\infty)]$ is divisible by $q$ for all but finitely 
many  $\ell\in \mathrm{Spl}(L)$. Furthermore $\ol{\eta}_\ell(\Gal(K(\mu_\infty))$ is a subgroup of $\ol{\eta}_\ell(\Gal(H'(\mu_\infty))$, and
hence $|\ol{\eta}_\ell(\Gal(K(\mu_\infty))|$ is divisible by $q$ for all $\ell\in \mathrm{Spl}(L)$. As $\ol{\eta}_\ell(\Gal(K(\mu_\infty))$
is a finite abelian group, it follows that $\Zz/q\in \FSQ(\ol{\eta}_\ell(\Gal(K(\mu_\infty)))$ for all $\ell\in \mathrm{Spl}(L)$.
Finally, by Chebotarev, the set $S:=\mathrm{Spl}(L)$ of prime numbers has a positive Dirichlet density.\hfill $\Box$

\begin{coro} \label{cmcoro}  
For every finite extension $E/K$ neither the restricted family $(\ol{\eta}_\ell|\Gal(E_{\cyc}))_{\ell\in\Ll}$ 
nor the restricted family $(\ol{\eta}_\ell|\Gal(E^\dagger_{\cyc}))_{\ell\in\Ll}$
is group theoretically independent. 
\end{coro}

{\em Proof.} Let $E/K$ be a finite extension. We apply Theorem \ref{cmdens} to $A_E/E$ in order to conclude
that $\Zz/2\in \FSQ(\ol{\eta}_\ell(\Gal(E_{\cyc}))$ for infinitely many $\ell\in\Ll$. As $\ol{\eta}_\ell(\Gal(E_{\cyc}))$ is a subgroup
of the finite abelian group $\ol{\eta}_\ell(\Gal(E_{\cyc}^\dagger))$, we conclude that also 
$\Zz/2\in \FSQ(\ol{\eta}_\ell(\Gal(E_{\cyc}^\dagger))$ for infinitely many $\ell\in\Ll$.\hfill $\Box$

\begin{rema} \label{cmrema} Corollary \ref{cmcoro} holds accordingly with $\ol{\eta}_\ell$ replaced by $\eta_\ell$. In fact, for every algebraic
extension $K'/K$ we have
$\FSQ(\ol{\eta}_\ell(\Gal(K')))\subset \FSQ({\eta}_\ell(\Gal(K')))$ for all $\ell\in\Ll$, because $A[\ell]=T_\ell(A)\otimes \Ff_\ell$ for
all $\ell\in\Ll$ and thus $\ol{\eta}_\ell(\Gal(K'))$ is a quotient of ${\eta}_\ell(\Gal(K'))$ for all $\ell\in\Ll$. 
\end{rema}
\end{appendices}

\subsection*{Acknowledgements} 
The Main Theorem of this paper agrees with a result  from the author's habilitation theses \cite{habil}. It continues the line of
investigations \cite{bible}, \cite{gp}, \cite{bgp}, \cite{cadoret} and is inspired largely by seminal work of 
Serre (cf. \cite{serre-cours1984}, \cite{serre-cours1985},  \cite{bible}). We acknowledge this with pleasure.
I want to thank Lior Bary-Soroker, Gebhard B\"ockle, Wojciech Gajda and
Cornelius Greither for interesting discussions. I acknowledge 
enlightning answers of Jim Humphreys and of an anonymous user to a question of mine
posed on the internet platform 
{\sc www.mathoverflow.net}. From this I learned the proof of Proposition \ref{gtprop}. I want to thank the
anonymous referee for a careful reading of the manuscript and for many helpful comments and suggestions.
Part of this work
was done during visits in Pozna\'n at Adam 
Mickiewicz University financed by NCN grant nr. UMO-2014/15/B/ST1/00128.

Sebastian Petersen\\
Fachbereich f\"ur Mathematik und Naturwissenschaften \\ 
Universit\"at Kassel\\
Wilhelmsh\"oher Allee 73\\
34121 Kassel, Germany\\
E-mail: petersen@mathematik.uni-kassel.de

\end{document}